\documentclass{amsart}

\usepackage{amscd}
\usepackage{amssymb, amsfonts}
\usepackage{amsrefs}
\input xy
\xyoption{all}

\newcommand{\Q}{{\mathbb Q}}
\newcommand{\C}{{\mathbb C}}
\newcommand{\Z}{{\mathbb Z}}

\newcommand{\gl}{\mathfrak{gl}}
\newcommand{\Gl}{\mathrm {GL}}
\newcommand{\ul}{{\mathfrak u}}
\newcommand{\ug}{{\mathrm U}}  
\newcommand{\fs}{{\mathfrak s}}

\newcommand{\oi}{{\bf \mathrm{O}}}

\newcommand{\ldp}{{\mathcal L}_{\mathrm {DP}}}
\newcommand\ldpo{{\mathcal L}_{\mathcal O}}
\newcommand\cmf[1]{{\mathcal C}(#1)}
\newcommand\ca{{\mathcal A}}
\newcommand\co{{\mathcal O}}
\newcommand\cb{{\mathcal B}}
\newcommand\ord{\mathrm{ord}}
\newcommand\ac{\overline{\mathrm{ac}}}
\newcommand\lef{\mathbb L}

\theoremstyle{plain}
\newtheorem{thm}{Theorem}
\newtheorem{theorem}[thm]{Theorem}

\theoremstyle{definition}
\newtheorem{remark}[thm]{Remark}

\title[]{Transfer to characteristic zero: appendix to ``Fundamental Lemma of Jacquet-Rallis in positive characteristics'' by Zhiwei Yun}
\author[]{Julia Gordon}

\begin{document} 

\maketitle

The purpose of this Appendix is to point out that the work of R. Cluckers, 
T. Hales and F. Loeser \cite{cluckers-hales-loeser} implies that 
Transfer principle of Cluckers and Loeser \cite{cluckers-loeser:fourier}
 applies to the version of the Fundamental 
Lemma proved in \cite{yun:FLJR}. Thus, the Conjectures 1.1.1 and 1.1.2 of 
\cite{yun:FLJR} 
are true when  $F$ is a local field of characteristic zero, with sufficiently 
large residue characteristic.

We need to emphasize that even though these Conjectures in the equal characteristic case are proved in \cite{yun:FLJR}
for the  fields $F$ of characteristic larger than $n$, 
the transfer principle leads to a slightly  weaker result for the fields of characteristic zero.  
Namely,  there exists (an algorithmically computable) constant $M$ such that  the Conjectures hold for 
the characteristic zero local fields  $F$ of 
residue characteristic larger than $M$. 
However, we hope that such a result 
is sufficient for some applications. 

Since this appendix is of an expository nature, the references often point not 
to the original sources, but to more expository articles.
All the references to specific sections, conjectures, 
and definitions that 
do not mention a source are to Z. Yun's article for which this appendix 
is written.

{\bf Acknowledgement.} This appendix emerged as a result of the AIM workshop on Relative Trace Formula and Periods of Automorphic 
Forms in August 2009. It is a pleasure to thank the organizers and participants of this workshop.  
I would like to emphasize that all the original ideas used and described 
here appear in the works of  R. Cluckers, T.C. Hales, and F. Loeser.  
I am grateful to R. Cluckers for a careful reading.

\section{Denef-Pas language}

The idea behind the approach to transfer described here  is to express 
everything involved in the statement of the Fundamental Lemma by means of 
formulas in a certain first-order language of logic (called the Denef-Pas 
language) $\ldp$ (see, e.g., \cite{cluckers-hales-loeser}*{Section 1.6} for the 
detailed definition), and then work with these formulas directly
 instead of the sets and functions described by them. 
Denef-Pas  language  is designed for valued fields. 
It is a {\it three-sorted language}, meaning that
it has three sorts of variables. Variables of the first sort run over 
the valued field, variables of the second 
sort run over the value group (for simplicity, we shall assume that  
the value group is $\Z$), and variables of the
third sort run over the residue field. 

Let us describe the set of symbols that, along with parentheses, 
the binary relation symbol `$=$' in every sort,  the standard logical 
symbols for conjunction, 
disjunction, and negation, and the quantifiers, are used to build formulas in 
Denef-Pas  language.

\begin{itemize}
\item 
In the valued field sort:  there are constant symbols `$0$' and `$1$', and  
 the symbols `$+$' and `$\times$' for the binary operations of 
addition and multiplication.
Additionally, there are symbols for two functions from the valued field sort:
`$\ord(\cdot)$' to denote a function from the valued field sort to the 
$\Z$-sort,  
and `$\ac(\cdot)$' to denote a function from
the valued field sort to the residue field sort.
These functions are called the {\it valuation map}, and the 
{\it angular component map}, respectively.  

\item In the residue field sort: there are constant symbols `$0$' and `$1$', 
and the binary operations symbols `$+$' and `$\times$' (thus, restricted to the residue field sort, this is the language of rings).

\item In the $\Z$-sort, there are `$0$' and `$1$', and the operation `$+$';  
additionally, for each $d=2,3,4,\dots$, there is 
a symbol `$\equiv_d$' to denote the binary 
relation $x\equiv y \text{ mod }d$. Finally, there is a binary relation symbol `$\ge$'.  
(This is Presburger language for the integers). 

\end{itemize}
Given a discretely valued field $K$ with a uniformizer of the valuation 
$\varpi$,  the functions `$\ord(\cdot)$' and `$\ac(\cdot)$' are interpreted as follows.
The function $\ord(x)$ stands for the valuation of $x$. 
It is  in order to provide the interpretation for the symbol `$\ac(x)$'
that a choice of the uniformizing parameter $\varpi$ (so that
$\ord(\varpi)=1$) is needed. 
If $x\in\co_{K}^{\ast}$ is a unit, there is a natural definition
of $\ac(x)$ -- it is the reduction of $x$ modulo the ideal $(\varpi)$.
For $x\neq 0$ in $K$, $\ac(x)$ is defined by $\ac(x)=\ac(\varpi^{-\ord(x)}x)$,
 and by definition, $\ac(0)=0$.

A formula $\varphi$ 
in $\ldp$ can be interpreted in any discretely valued field, once a 
uniformizer of the valuation is chosen, in the sense that given a valued field 
$K$ with a uniformizer $\varpi$ and the residue field $k_K$, 
one can allow the free variables of $\varphi$ to range over $K$, $k_K$, 
and $\Z$, respectively, according to their sort (naturally, the variables 
bound by a quantifier then also range over $K$, $k_K$, 
and $\Z$, respectively). 
Thus, any  discretely valued field 
is a structure for Denef-Pas language.

\section{Constructible motivic functions}
In the foundational papers  \cite{cluckers-loeser}, 
\cite{cluckers-loeser:fourier}, R. Cluckers and F. Loeser developed the theory of motivic integration for functions defined by means of formulas in Denef-Pas 
language, and proved a very general Transfer Principle.
We refer to \cite{cluckers-hales-loeser} and \cite{cluckers-loeser:expo}
for the introduction to this subject and all definitions. Note
that the article \cite{cluckers-hales-loeser} is self-contained and 
essentially covers everything in this appendix.

Here we need to use the terms ``definable subassignment'' and 
``constructible motivic function''. 
Let $h[m,n,r]$ be the functor from the category of fields to the category of 
sets defined by 
\begin{equation*}
h[m,n, r](K)=K((t))^m\times K^n \times \Z^r.
\end{equation*}
The term {\sl subassignment} was first introduced in 
\cite{denef-loeser:p-adic}. Given a functor $F$ from some category $C$ 
to $\underline{\text{Sets}}$,  a subassignment $X$ of $F$ is a collection of subsets 
$X(A)\subset F(A)$ for every object $A$ of $C$. 
A definable set is a set that can be described by a formula in Denef-Pas 
language, and a subassignment $X$ of the functor $h[m,n,r]$ is called definable
if there exists a formula $\varphi$ in Denef-Pas language with $m$ free variables of the valued field sort, $n$ free variables of the residue field sort, 
and $r$ free variables of the value group sort, 
such that for every field $K$, the set $X(K)$ is exactly the set of points
in $K((t))^m\times K^n \times \Z^r$ where $\varphi$ takes the value `true'.
Note that there are slightly different variants of Denef-Pas language, depending on the sets of coefficients for a formula $\varphi$  
allowed in every sort 
(the smallest set of coefficients is $\Z$ in every sort; however, one can add constant symbols that can later be used as coefficients -- one such variant 
will be discussed below). 
We emphasize, 
however, that regardless of the variant, the coefficients come from a 
fixed set, and are 
independent of $K$.  
Definable subassignments form a Boolean algebra in a natural way, and this 
algebra is the replacement, in the theory of motivic integration, for the 
Boolean algebra 
of measurable sets in the traditional measure theory.

For a definable subassignment $X$, the ring  
of the so-called {\sl constructible motivic  
functions} on $X$, denoted by $\cmf{X}$, is defined in 
\cite{cluckers-loeser}.  The elements of $\cmf{X}$ are, 
essentially, formal constructions defined using the language $\ldp$. 
The main feature of constructible 
motivic functions is specialization to functions on discretely valued fields. 
Namely, let $f\in \cmf{X}$. 
Let $F$ be a non-Archimedean local field 
(either of characteristic zero or of positive characteristic). 
Let $\varpi$ be the uniformizer of the 
valuation on $F$. Given these data, one gets a specialization 
$X_F$ of the subassignment $X$ to $F$, which is a definable subset of 
$F^m\times k_F^n\times \Z^r$ for some $m,n,r$, and 
the constructible motivic function $f$ specializes to a 
$\Q$-valued function $f_F$ on $X_F$, for all fields
$F$ of residue characteristic bigger than a constant that depends only on the 
$\ldp$-formulas defining $f$ and $X$.  
As explained in \cite{cluckers-hales-loeser}*{Section 2.9}, one can tensor the 
ring $\cmf{X}$ with $\C$, and then the specializations $f_F$ of 
elements of $\cmf{X}\otimes \C$ form a $\C$-algebra of functions on $X_F$. 

\section{Integration and Transfer Principle} 
In \cite{cluckers-loeser}, Cluckers and Loeser defined a class $\mathrm I C(X)$ of {\sl integrable} 
constructible motivic functions, closed under  integration with respect to 
parameters (where integration is with 
respect to the {\sl motivic measure}). Given a local field $F$ with a choice 
of the uniformizer, these functions specialize to 
integrable (in the classical sense) functions on $X_F$, and motivic 
integration specializes to the 
classical integration with respect to an appropriate Haar measure, when the 
residue 
characteristic of $F$ is sufficiently large.

From now on, we will use the variant of the theory of motivic integration 
with coefficients in the ring of integers
of a given global field.
Let $\Omega$ be a global field with the ring of integers ${\mathcal O}$. 
Following \cite{cluckers-loeser:expo}, we denote by 
$\ca_{\co}$ the collection of all $p$-adic completions of all finite 
extensions of $\Omega$, and by $\cb_{\co}$ the set of all local fields of 
positive characteristic that are $\co$-algebras.
Let $\ca_{\co, M}$ (resp.,  $\cb_{\co, M}$) be the set of all local fields
$F$ in $\ca_{\co}$ (resp., $\cb_{\co}$) such that the residue field $k_F$ has 
characteristic larger than $M$. 
Let $\ldpo$ be the variant of Denef-Pas language with coefficients in 
${\mathcal O}[[t]]$ (see \cite{cluckers-loeser:expo}*{Section 6.7} for the 
precise definition). This means, roughly, 
that a constant symbol for every element of $\mathcal O[[t]]$ is added to the valued field sort, so that a formula in $\ldpo$ is allowed to have coefficients in ${\mathcal O}[[t]]$ in the valued field sort, coefficients in $\Omega$ in the residue field sort, and coefficients in $\Z$ 
in the value group sort.
  
Then the Transfer Principle can be stated as follows. 

\begin{theorem}(Abstract Transfer Principle, cf. 
\cite{cluckers-hales-loeser}*{Theorem 2.7.2}).
Let $X$ be a definable subassignment, and let $\varphi$ be a constructible
(with respect to the language $\ldpo$) motivic function on $X$. 
Then there exists $M>0$ such that for every 
$K_1, K_2 \in \ca_{\co, M}\cup \cb_{\co, M}$ with $k_{K_1}\simeq k_{K_2}$,
$$
\varphi_{K_1}=0 \quad \text{if and only if}\quad \varphi_{K_2}=0.
$$
\end{theorem}

\begin{remark}
In fact, the transfer principle is proved in 
\cite{cluckers-loeser:fourier} for an even richer class of functions, 
called constructible motivic exponential functions, that 
contains additive characters of the field along with the constructible motivic functions. However, we do not discuss it here since the characters are not 
needed in the present setting. 
\end{remark}

The goal of this appendix is to check that the Conjectures proved in 
this article can be expressed as equalities between specializations 
of constructible motivic functions. We emphasize that all the 
required work 
is actually done 
in \cite{cluckers-hales-loeser}, here we just check that it indeed applies 
in the present situation.

\section{Definability of all the ingredients}
Here we go through Section 2.1 and check that every object 
appearing in it is definable.

\subsection{The degree two algebra $E/F$}\label{subsub:E} 
Following \cite{cluckers-hales-loeser}*{Section 4}, we fix, once and for all,
a $\Q$-vector space $V$ of dimension $n$, and fix a basis 
$e_0,\dots, e_n$ of $V$ over $\Q$. 

As in \cite{cluckers-hales-loeser}*{Section 3.2}, we introduce a parameter 
(which we denote by $\epsilon$) that will appear in all the formulas that 
involve an unramified quadratic extension of the base field. 
We think of $\epsilon$ as a non-square unit, and denote by $\Lambda$ be the 
subassignment of $h[1,0,0]$ defined by the formula 
`$\ord(\epsilon)=0 \wedge \nexists x: x^2=\epsilon$'. 
From now on, we only consider the relative situation: all the subassignments
we consider will come with a fixed projection morphism to $\Lambda$
(in short, we are considering the  category of definable subassignments 
over $\Lambda$, see \cite{cluckers-loeser}*{Section 2.1}). 
That is, we replace all the constructions that depend on 
an unramified quadratic extension $E/F$ (such as the unitary group), with the 
family of isomorphic objects parameterised by a non-square unit $\epsilon$ in 
$F$. 
Now, imagine that we fixed the basis $(1, \sqrt{\epsilon})$ for the quadratic 
extension $E$. Then $E$ can be identified with $F^2$ via this basis; so from 
now on we shall think of the elements of $E$ as pairs of variables that 
range over $F$. The nontrivial Galois automorphism  $\sigma$ 
of $E$ over $F$ now can be 
expressed as a $2\times 2$-matrix with entries in $F$, and can be used in the 
expressions in Denef-Pas language.

The nontrivial quadratic character $\eta_{E/F}$ can be expressed by a 
Denef-Pas formula `$\eta_{E/F}(x)=1 \Leftrightarrow \exists (a,b)\in F^2 : (a^2+\epsilon b^2=x)$', or, simply, `$\eta_{E/F}(x)=1 \Leftrightarrow\ord(x)\equiv 0 \mod 2$'. 

In the case $E/F$ split, we just treat elements of $E$ as pairs of elements of 
$F$. 

\subsection{The groups and their Lie algebras}\label{subsub:lie}
In Sections 2.1, 2.2, one starts out with free 
${\mathcal O}_F$-modules $W$ and $V$, and then proceeds to choose a basis 
vector $e_0$ with certain properties. We shall reverse the thinking here: 
we fix a basis $e_0, e_1, \dots e_n$, and fix the the dual basis $e_0^{\ast},\dots, e_n^{\ast}$, such that $e_0^{\ast}(e_0)=1$, and  
such that the Hermitian form $(\cdot ,\cdot )$ on $V$, with respect to this basis, corresponds to a matrix with entries in the set $\{0, \pm 1\}$, and $(e_0, e_0)=1$, and 
$(e_0, e_j)=0$ for $1\le j\le n$. 
We let $W$ be the span of the vectors $e_1, \dots, e_{n-1}$. 
With this choice of basis, we  think of the elements of $\gl_n$ 
as $n^2$-tuples of variables $A=(a_{ij})$.
(Formally speaking, we identify $\gl_n$ with the definable subassignment
$h[n^2,0,0]$.)
All the split algebraic groups are, naturally, defined by polynomial 
equations in these variables, and thus can be replaced with definable 
subassignments 
of $h[n^2,0,0]$. The embedding 
$\Gl_{n-1}\hookrightarrow\Gl_n$ where 
$$
A\mapsto \left(\begin{matrix} A&{}\\{}& 1 \end{matrix}\right)
$$
is, clearly, definable.

To find the definable subassignments that specialize
to ${\mathfrak s}_n$, ${\mathfrak u}_n$, and $\ug_n$, we introduce the parameter
$\epsilon$ as above in Section \ref{subsub:E}.
Then $\fs_n$ naturally becomes a definable subassignment of 
$h[2n^2, 0,0]\times \Lambda\subset h[2n^2+1, 0,0]$. 
Indeed, as discussed above, the Galois automorphism $\sigma$ can be used in 
$\ldp$-expressions when we think of the elements of $E$ as pairs of 
$F$-variables: we replace each variable $a_{ij}$  ranging over $E$ 
with a pair of variables
$(x_{ij}, y_{ij})$ ranging over $F$.
The Hermitian form that is used to define the unitary group, given the 
choice of the basis, gives rise to polynomial equations in 
$(x_{ij}, y_{ij})$ that define the unitary group.  Hence, $\ul_n$ and $\ug_n$ can also be replaced 
with definable subassignments of $h[2n^2+1, 0,0]$. 

\subsection{The invariants}\label{subsub:inv}
By definition, $a_i(A)$ are the coefficients of the characteristic 
polynomial of $A$, and in particular, they are polynomial expressions in the matrix entries of $A$, and therefore they are given by terms in $\ldp$, and the map
$A\mapsto (a_i(A))_{0\le i \le {n-1}}$ is definable (recall that a function is called definable if its graph is a definable set).

First, let us consider the case when $E/F$ is a field extension.

The linear functional $e_0^{\ast}$ on $V$ (defined in Section 2.1), with our choice of the bases, is just the covector $(1, 0, \dots, 0)$. 
Then the invariants $b_i(A)$ of Section 2.2 
are also given by terms in $\ldp$. 

The vectors $A^ie_0$ are, clearly, just columns of polynomial expressions in
the matrix entries $(x_{ij}, y_{ij})$ of $A$.  
The condition that a collection of vectors forms a basis of a given 
vector space is  a predicate in $\ldp$. 
Hence, the set of semisimple elements in $\gl_n(E)$ that are strongly regular
with respect to $\Gl_{n-1}(E)$-action (in the sense of Definition 2.2.1)
is a specialization (to $F$) of a definable subassignment of $h[2n^2+1, 0,0]$.

We observe that $\Delta_{a,b}=\det(e_0^{\ast} A^{i+j}e_0)_{0\le i,j\le n}$ 
(of Definition 2.2.3) is also 
a polynomial expression in $(x_{ij}, y_{ij})$. 

Recall the subassignment $\Lambda$ from the previous subsection 
that specializes to the domain for a parameter $\epsilon$ defining 
the extension $E$.
Since the image of a definable subassignment under a definable morphism 
is a definable subassignment, we have the definable subassignment 
${\mathcal P}$ over $\Lambda$, which we will 
denote by ${\mathcal P} \to \Lambda$, that corresponds to the set of pairs 
$(a,b) \in E^{2n}$ that are invariants of some strongly regular element of 
$\gl_n(E)$.  Mote precisely, ${\mathcal P}$ is a subassignment of 
$\Lambda\times h[4n,0,0]$ that satisfies the condition 
that there exists $N>0$, such that for every local 
field $F\in \ca_{\co, N}\cup \cb_{\co, N}$, 
for every $\epsilon \in \Lambda_F$, the fibre 
${\mathcal P}_{\epsilon}$ of ${\mathcal P}$ at $\epsilon$ specializes to the 
set of pairs $(a,b)$ that are invariants of some $A\in \gl_n(E)$,  
strongly regular
with respect to $\Gl_{n-1}(E)$-action (in the sense of Definition 2.2.1),
where 
$E$ is the field extension corresponding to $\epsilon$.  

If $E/F$ is split, the same argument works, except there is no need to 
consider the relative situation over $\Lambda$. 

Since we have a symbol for the $F$-valuation in $\ldp$, the parameter 
$\nu(A)$ of Definition 2.2.2 is also an expression in $\ldp$.

\subsection{The orbital integrals}
Since the quadratic character $\eta_{E/F}$ only takes the values $\pm 1$, 
we can break the orbital integral 
${\oi}_A^{\Gl_{n-1,\eta}}\left({\bf 1}_{\fs_n(\co_F)}\right)$ into 
the difference of two integrals:
\begin{equation*}
\begin{aligned}
{\oi}_A^{\Gl_{n-1,\eta}}\left({\bf 1}_{\fs_n(\co_F)}\right)& =
\int_{\Gl_{n-1}(F)\cap \{g\mid \eta_{E/F}(\det g)=1\}}{\bf 1}_{\fs_n(\co_F)}(g^{-1}Ag)\,dg\\
& -\int_{\Gl_{n-1}(F)\cap \{g\mid \eta_{E/F}(\det g)=-1\}}
    {\bf 1}_{\fs_n(\co_F)}(g^{-1}Ag)\,dg.
\end{aligned}
\end{equation*}

By the remarks in Section \ref{subsub:E} above, both domains of integration are 
definable sets. 
For each point $A$ in the subassignment of strongly regular elements,
${\bf 1}_{\fs_n(\co_F)}(g^{-1}Ag)$ is, by Section \ref{subsub:lie} above, a specialization of a constructible motivic function of $g$.
We need to briefly discuss the normalization of the measures.
The $p$-adic measure to which the motivic measure specializes is the so-called 
Serre-Oesterl\'e measure, defined in \cite{oesterle}. 
 Serre-Oesterl\'e measure on 
a classical group $G$ is the Haar measure such that the volume of the maximal compact subgroup 
is  $q^{\dim G}$.
Hence, the Haar measure $dg$ differs from Serre-Oesterl\'e measure by a 
factor of $q^{-(n-1)^2}$, where $q$ is the cardinality of the residue field, since, as in \cite{jacquet-rallis}, the Haar measures here 
are chosen so that the standard maximal compact subgroups have volume $1$. 
This factor is the specialization of the (constant) 
constructible motivic function $\lef^{-(n-1)^2}$ 
(see e.g.\cite{cluckers-hales-loeser}*{Section 2.3} 
for the discussion of the symbol 
$\lef$). 
We conclude that 
${\oi}_A^{\Gl_{n-1,\eta}}\left({\bf 1}_{\fs_n(\co_F)}\right)$ is a specialization of a constructible motivic function of $A$. 

By a similar inspection, we see that the integral 
$\oi_{A'}^{U_{n-1}}\left({\bf 1}_{\ul_n(\co_F)}\right)$ is a specialization of a 
constructible motivic function 
of $A'$, and thus so is 
the right-hand side of Conjecture 1.1.1 (1). 

Finally, recall the subassignment ${\mathcal P}$ from 
Section \ref{subsub:inv} that specializes to the set of invariants. 
Consider the subassignment ${\mathcal X}$ of $\fs_n\times \ul_n $ defined by:
$(A, A')\in {\mathcal X}$ if and only if $A$ and $A'$ have the same invariants.
Since as we discussed above, the map that maps $A$ to its collection of
invariants is a definable map, this is a definable subassignment (note that it has a natural projection to ${\mathcal P}$). We have shown that the difference of the left-hand side and the right-hand side of Equation (1)
in Conjecture 1.1.1 is a constructible motivic 
function on ${\mathcal X}$. Therefore, the Transfer Principle applies to it.

By inspection, all the ingredients of all the other variants of 
Conjecture 1.1.1 and Conjecture 1.1.2 are definable in the language $\ldpo$, and hence the Transfer Principle applies in all these cases.

\end{document}